\documentclass[12pt]{article}
\usepackage{cite,amsfonts,amsmath,subeqn}

\def\be{\begin{equation}}
\def\eqn#1{\be\label{#1}}

\def\bea{\begin{eqnarray}}

\def\eea{\end{eqnarray}}

\def\nn{\nonumber}

\newcommand{\eqna}[1]{\begin{subequations} \label{#1}
\begin{eqnarray}}
\def\eena{\end{eqnarray}
\end{subequations}}

\def\dia{$\ \diamondsuit$} 
\def\PR{{\it Proof: }}
\def\np{\vfil\eject}
\def\nl{\hfil\break}
\def\nt{\noindent}
\def\id{{\bf 1}}
\def\ta{\tilde a}\def\td{\tilde d}
\def\tA{\tilde A}\def\tD{\tilde D}
\def\hR{\hat R} 
\def\g{\gamma} \def\z{\zeta} 
\def\d{\delta}
\def\ve{\varepsilon}
\def\cU{{\cal U}} \def\un{\,{{\cal U}}_1\,}
\def\ut{\,{{\cal U}}_3\,} \def\ur{\,{{\cal U}}_2\,}
\def\cA{{{\cal A}}} 
\def\an{{\,{\cal A}}_1\,}
\def\at{{\,{\cal A}}_3\,}
\def\ar{{\,{\cal A}}_2\,}
\def\cD{{\cal D}} \def\cP{{\cal P}} \def\cQ{{\cal Q}}
\def\lg{\left\langle}
\def\rg{\right\rangle}
\def\rra{\longrightarrow}
\def\ra{\rightarrow}
\def\lra{\leftrightarrow}
\def\ha{\textstyle{\frac{1}{2}}} 
\def\mt{\mapsto}
\newcommand{\CC}{\mbox{${\mathbb C}$}}

\newcommand{\ZZ}{\mbox{${\mathbb Z}$}}
\newcommand{\NN}{\mbox{${\mathbb N}$}}

\newcommand{\Rtriang}[9]{
 {\left(
 \begin{array}{cccc}
 1 & \ #1\ & \ #2\ & \ #3 \ \cr
 & #4 & #5 & \ #6\ \cr
 & & \ #7\ & \ #8\ \cr
 & & & #9
 \end{array}
 \right)
 }}

\normalbaselines 
\begin{document}
\pagestyle{empty}

\begin{center}

 \textsf{\LARGE Duality for Exotic Bialgebras}

\vspace{7mm}

{\large D.~Arnaudon$^{a,}$\footnote{Daniel.Arnaudon@lapp.in2p3.fr}, 
 ~A.~Chakrabarti$^{b,}$\footnote{chakra@cpht.polytechnique.fr},\\[2mm] 
 V.K.~Dobrev$^{c,d,}$\footnote{vladimir.dobrev@unn.ac.uk,dobrev@inrne.bas.bg} 
~and~ S.G.~Mihov$^{d,}$\footnote{smikhov@inrne.bas.bg}}

\vspace{5mm}

\emph{$^a$ Laboratoire d'Annecy-le-Vieux de Physique Th{\'e}orique LAPTH}
\\
\emph{CNRS, UMR 5108, associ{\'e}e {\`a} l'Universit{\'e} de Savoie}
\\
\emph{LAPTH, BP 110, F-74941 Annecy-le-Vieux Cedex, France}
\\
\vspace{3mm}
\emph{$^b$ Centre de Physique Th{\'e}orique, CNRS UMR 7644}
\\
\emph{Ecole Polytechnique, 91128 Palaiseau Cedex, France.}
\\
\vspace{3mm}
\emph{$^c$ School of Computing and Mathematics}
\\
\emph{University of Northumbria}
\\
\emph{Ellison Place, Newcastle upon Tyne, NE1 8ST, UK}\\ 
\vspace{3mm} 
\emph{$^d$ Institute of Nuclear Research and Nuclear Energy} 
\\
\emph{Bulgarian Academy of Sciences}
\\
\emph{72 Tsarigradsko Chaussee, 1784 Sofia, Bulgaria}
\footnote{permanent address for V.K.D.}
\\
\vspace{3mm}

\end{center}

\vspace{.8 cm}

\begin{abstract}
In the classification of Hietarinta, three triangular $4\times 4$
$R$-matrices lead, via the FRT formalism, to matrix bialgebras
which are not deformations of the trivial one. In this paper, we
find the bialgebras which are in duality with these three exotic
matrix bialgebras. We note that the $L-T$ duality of FRT is not
sufficient for the construction of the bialgebras in duality. We
find also the quantum planes corresponding to these bialgebras
both by the Wess-Zumino R-matrix method and by Manin's method.
\end{abstract}

\vfill

\rightline{LAPTH-823/00}
\rightline{CPHT-S 003.0101}
\rightline{UNN-SCM-M-00-13}
\rightline{INRNE-TH-00-06}
\rightline{math.QA/0101160}
\rightline{December 2000}

\newpage
\pagestyle{plain}
\setcounter{page}{1}

\section{Introduction}
\label{sect:intro}
\setcounter{equation}{0}

Until recently it was not clear how many distinct quantum group
deformations are admissible for the group $\ GL(2)\ $ and the
supergroup $\ GL(1|1)\ $. For the group $\ GL(2)\ $ there were
the well-known standard $GL_{pq}(2)$ \cite{DMMZ} 
and nonstandard (Jordanian) $GL_{gh}(2)$ \cite{Ag} two-parameter
deformations. For the supergroup $\
GL(1|1)\ $ there were the standard $GL_{pq}(1|1)$
\cite{HiRi,DaWa,BuTo} and the hybrid (standard-nonstandard)
$GL_{qh}(1|1)$ \cite{FHR} two-parameter deformations. 
(Various aspects of these matrix quantum (super-)group deformations 
were studied in, e.g., \cite{Manin,CFFS,SWZ,EOW,OW,Sudbery,Zak,Dod,Kup,Ohn,%
FaZa,Vla,Ogiev,Kar,SchSch,BCHOS,JagJeu,ABCD,Dobh,ChaJag,ADM,KP,ACCt,JacCor,%
ACCS,CQ,BHP,Aiz,CC,AACDM,Chak,LMO}.) Recently, in \cite{AACDM} it was
shown that the list of these four deformations is exhaustive
(refuting a long standing claim of \cite{Kup} (supported also in
\cite{BHP,LMO}) for the existence 
of a hybrid (standard-nonstandard) two-parameter deformation of
$GL(2)$). In particular, it was shown that these four
deformations match the distinct triangular $4\times 4$
$R$-matrices from the classification of \cite{Hietarinta} 
which are deformations of the trivial $R$-matrix (corresponding
to undeformed $GL(2)$). 

The matching mentioned above was done 
 by applying the FRT formalism \cite{FRT} to these $R$-matrices.
While applying this we noticed that a particular $R$-matrix, 
namely, the one denoted by $\ R_{H2,3}\ $ in \cite{Hietarinta}, 
gives different results depending on the range of the three
parameters $\ h_1\,,\,h_2\,,\,h_3\ $ it depends of. Only one of
the ranges, namely, $\ h_1 = -h_2= h \ $, $\ h_3 = -h^2\ $ contains the 
zero point which gives the trivial $R$ matrix. (It is a partial 
case of the two-parameter $\ g,h\ $ Jordanian deformation for $\
g=-h\ $.) The other two ranges are given by $\ h_1 = -h_2= h \ $,
$\ h_3 \neq -h^2\ $ $(*)$ and $\ h_1 \neq -h_2$ $(**)$. Thus, the
$R$-matrices obtained while varying through these ranges are not
deformations of the trivial $R$ matrix, and also are distinct
between each other. This analysis revealed altogether three distinctly
different triangular $R$-matrices which are ~{\it not}~
deformations of the trivial $R$-matrix. 
In this way, in \cite{AACDM} were obtained three new matrix bialgebras
which are ~{\it not}~ deformations of the classical algebra 
of functions over the group $\ GL(2)$ or the supergroup
$GL(1|1)$. These new matrix bialgebras, which we now call ~{\it exotic}~ are
very interesting and deserve 
further study. One of the first problems when dealing with such
matrix bialgebras is to find the bialgebras with which they are in
duality, since some of the structural 
characteristics are more transparent for the duals. The
bialgebras in duality are also the
interesting objects with respect to the development of the
representation theory. 

This is the problem we solve in this paper. We find the
bialgebras which are in duality with the three exotic matrix
bialgebras found in \cite{AACDM}. We then find the quantum planes
corresponding to these bialgebras by the Wess-Zumino
R-matrix method \cite{WZ} (cf. also \cite{Schwenk}). 
For the latter we find the minimal polynomials\ pol$(\cdot)$ in
one variable such that \ pol$(\hR)=0$ is the lowest order
polynomial identity satisfied by the singly permuted $R$-matrix
$\hR \equiv P\,R$ ($P$ is the permutation matrix).  These
minimal polynomials indeed separate the three cases of $\
R_{H2,3}\ $ mentioned above. Namely, in case $(*)$ we find a
cubic minimal polynomial, while in the case $(**)$ it is quartic,
(cf.  (\ref{hwrm}),(\ref{mpol})). (Recall that the corresponding minimal
polynomial in the Jordanian case is only quadratic.) We find also
the quantum planes by Manin's method \cite{Manin}.

The paper is organized as follows. In Section 2 we give a the
overall general setting. Sections 3,4,5, are devoted to the
separate study of the three exotic bialgebras. In each case we
first give a more detailed (than in \cite{AACDM}) 
picture of the structure of these
bialgebras, then we construct the bialgebra in duality, noting
the bearing it has on the initial bialgebra; finally we show
consistency of this approach to duality with the FRT one, noting
the failings of the latter for these exotic bialgebras. In
Section 6 we construct the quantum planes corresponding to the
three matrix bialgebras. Section 7 
contains conclusions and outlook.

\section{Exotic bialgebras: general setting} 
\setcounter{equation}{0}

In this paper we consider the three exotic matrix bialgebras (obtained in
\cite{AACDM}) which are ~{\it not}~ deformations of the classical algebra 
of functions over the group $\ GL(2)$ or the supergroup
$GL(1|1)$. In all three cases these are 
unital associative algebras generated by four elements
$\ a,b,c,d$ which are not deformations of the classical algebra 
of functions over the group $\ GL(2)$ (or over the supergroup
$GL(1|1)$). This is evident also from 
the algebraic relations which we give separately below. 
The coalgebraic relations are the classical ones:
\eqn{coop}
\d\, (T) \ =\ T \otimes T \ , \qquad 
\ve\, (T) \ =\ \id_2\ \equiv \ \left(\begin{array}{cc}
1 & 0 \cr 0 & 1
\end{array} \right)
 \ , \qquad 
T\ = \ \left( \begin{array}{cc}
a & b \cr c & d 
\end{array} \right)
\end{equation} 
or explicitly:
\eqn{copr}
\d\ \left( \begin{array}{cc}
a & b \cr c & d 
\end{array} \right)
 \ =\ \left( \begin{array}{cc}
 a \otimes a + b \otimes c & a \otimes b + b \otimes d \cr 
c \otimes a + d \otimes c & c \otimes b + d \otimes d
\end{array} \right) 
\end{equation} 
\eqn{cou} 
\ve\ \left( \begin{array}{cc}
a & b \cr c & d 
\end{array} \right)
 \ =\ \left(\begin{array}{cc}
1 & 0 \cr 0 & 1
\end{array} \right)
 \end{equation}
However, the bialgebras under consideration are not Hopf
algebras, as we shall show in each particular case. 

\section{Exotic bialgebras: case 1} 
\setcounter{equation}{0}

\subsection{Bialgebra relations}

In this Section we consider the matrix bialgebra, denoted here by
$\ \cA_1\ $, which is obtained by applying the RTT relations of 
\cite{FRT}: 
\eqn{rtt} R\ T_1\ T_2 \ \ =\ \ T_2\ T_1\ R \ \ , \end{equation} 
where \ $T_1 \ =\ T\, \otimes\, \id_2$\ , \ $T_2 \ =\ \id_2 \,
\otimes\, T$, for the case when $\ R\ =\ R_1\ $: 
\eqn{rhu} 
R_1\ =\ \Rtriang{h}{-h}{h_3}
 {1}{0}{-h}
 {1}{h}
 {1}
 \quad , \qquad h_3\neq -h^2 
\end{equation} 
This $R$-matrix together with the condition on the parameters is
one of the special cases (mentioned in the Introduction) of the
$R$-matrix denoted by $\ R_{H2,3}\ $ in
\cite{Hietarinta}. 
The algebraic relations of $\ \an\ $ obtained in this way are
given by formulae (5.11) of \cite{AACDM}, namely: 
\begin{eqnarray}
 \label{eq:RTTh23b}
 &&c^2 \ =\ 0, \qquad ca \ =\ ac\ =\ 0, 
 \qquad dc \ =\ cd\ =\ 0, \nonumber\\
 && da \ =\ ad, \qquad cb \ =\ bc, \qquad a^2 \ =\ d^2 \nonumber\\
 && ab \ =\ ba + h(a^2 +bc-ad), \qquad
 db \ =\ bd - h(a^2 +bc-ad). 
\end{eqnarray}
Note that the constant $h_3$ does not enter the above 
relations. 

Note that this bialgebra is not a Hopf algebra. Indeed, suppose
that it is and there is an antipode $\ \gamma\ $, then we use one
of the Hopf algebra axioms:
\eqn{haxi} m \circ ({\rm id} \otimes \g) \circ \d \ =\ i \circ \ve 
\end{equation} 
as maps $\cA \ra \cA$, 
where $m$ is the usual product in the algebra: $m(Y\otimes Z)
 = YZ , Y,Z\in\cA$ and $i$ is the natural embedding of the
number field $F$ into $\cA$ : $i(c) = \mu 1_\cA , \mu\in
F$. Applying this to the element $\ d\ $ we would have:
$$c\,\g (b) + d\,\g(d) = 1_\cA $$ 
which leads to contradiction after multiplying from the left by
$\ c\ $ (one would get $0=c$).

The algebra $\an$ has the following PBW basis: 
\begin{equation}
 \label{eq:PBWb}
b^n a^k d^\ell \,, \quad b^n c\,, \qquad 
n,k\in \ZZ_+ \,, \quad \ell=0,1.
\end{equation}

The last line of (\ref{eq:RTTh23b}) strongly suggests the
substitution: 
\eqn{subs} \tilde{a} \ =\ \ha (a+d), \ \ \ \ \tilde{d} \ =\ \ha 
(a-d)\ , \end{equation} \nt 
so that the new algebraic relations and PBW basis are: 
\begin{eqnarray}\label{nrel}
&&c^2 \ =\ 0, \ \ \ \ \tilde{a}c \ =\ c\tilde{a} \ =\
\tilde{d}c \ =\ c\tilde{d} \ = \ \tilde{a}\tilde{d} \ =\
\tilde{d}\tilde{a} \ =\ 0 ,\ \ \ \ cb \ =\ bc, \nn\\ 
&&
\tilde{a}b \ =\ b\tilde{a}, \ \ \ \ \tilde{d}b \ =\ b\tilde{d} +
2h\tilde{d}^2 + hbc \eea 
\begin{equation}
 \label{eq:PBWbb}
b^n \ta^k\,, \quad b^n \td^\ell \,, \quad b^n c\,, \qquad 
n,k\in \ZZ_+ \,, \quad \ell\in\NN .
\end{equation}

The coalgebra relations become: 
\bea 
&&\d \left( \begin{array}{c}
\ta \cr\cr b \cr\cr c \cr\cr \td 
\end{array} \right) \ =\ 
\left( \begin{array}{c}
\ta \otimes \ta + \td \otimes \td + \ha\, b\otimes c + \ha\, c\otimes
b \cr \cr 
\ta \otimes b + \td \otimes b + b \otimes \ta - b \otimes \td
\cr \cr 
c \otimes \ta + c \otimes \td + \ta \otimes c - \td \otimes c
\cr\cr 
 \ta \otimes \td + \td \otimes \ta + \ha\, b\otimes c - \ha\, c\otimes b 
\end{array} \right) \label{capr} \qquad \\ \nn\\ 
&&\ve \left( \begin{array}{cc} 
\ta & b \cr c &\td \end{array} \right) \ =\ 
\left( \begin{array}{cc} 1 & 0 \cr 0 & 0 
\end{array} \right) \label{cau} \eea

\subsection{Duality}

Two bialgebras \ $\cU , \cA$\ are said to be \ {\it in duality}\ 
\cite{Abe} if there exists a doubly nondegenerate bilinear form 
\eqn{dua} \lg \ , \ \rg \ :\ \cU \times \cA \rra \CC \ , \ \ \ 
\lg \ , \ \rg \ :\ (u,a) \ \mt \ \lg \ u \ , \ a \ \rg
\ , \ u \in\cU \ , \ a \in \cA \end{equation} 
such that, for \ $u,v\in\cU\ , a,b\in\cA$\ :
\eqna{dub} 
&\lg \ u \ , \ ab \ \rg \ =\ \lg \ \d_\cU(u) \ ,
\ a \otimes b \ \rg \ , \ \ \ \lg \ uv \ , \ a \ \rg \ =\ 
\lg \ u \otimes v \ , \ \d_\cA (a) \ \rg \qquad \\ 
&\lg 1_\cU , a \rg \ =\ \ve_\cA (a) \ , \ \ \ \lg u ,
1_\cA \rg \ =\ \ve_\cU (u) \eena \nt 
Two Hopf algebras \ $\cU , \cA$\ are said to be \ {\it in duality}\ 
\cite{Abe} if they are in duality as bialgebras and if 
\be \lg \g_\cU (u) , a \rg \ =\ \lg u , \g_\cA (a)
\rg \end{equation} 

It is enough to define the pairing (\ref{dua}) between the generating
elements of the two algebras. The pairing between any other
elements of \ $\cU, \ \cA$\ follows then from relations (\ref{dub}) and
the standard bilinear form inherited by the tensor product. 

The duality between two bialgebras or Hopf algebras may be used
also to obtain the unknown dual of a known algebra. For that it
is enough to give the pairing between the generating elements of
the unknown algebra with arbitrary elements of the PBW basis of
the known algebra. Using these initial pairings and the duality
properties one may find the unknown algebra. One such possibility
is given in \cite{FRT}. However, their approach is
not universal. In particular, it is not enough for the algebras
considered here, (as will become clear) and will be used only as
consistency check. 

Another approach was initiated by Sudbery \cite{Sudbery}. He
obtained \ $U_q(sl(2)) \otimes U(u(1))$\ as the algebra of
tangent vectors at the identity of \ $GL_q(2)$. The initial pairings
were defined through the tangent vectors at
the identity. However, such calculations become very difficult
for more complicated algebras. Thus, in \cite{Dod} a
generalization was proposed in which the initial pairings are postulated to
be equal to the classical undeformed results. This generalized
method was applied in \cite{Dod} to the standard two-parameter
deformation $GL_{p,q}(2)$, (where also Sudbery's method was
used), then in \cite{DPa} to the multiparameter deformation of
$GL(n)$, in \cite{DPb} to the matrix quantum Lorentz group of
\cite{QL}, in \cite{ADM} to the Jordanian two-parameter
deformation $GL_{g,h}(2)$, in \cite{FHR} to the hybrid two-parameter
deformation of the superalgebra $GL_{q,h}(1|1)$, 
in \cite{DT} to the multiparameter
deformation of the superalgebra $GL(m/n)$. 
(We note that the dual of $GL_{p,q}(2)$ was obtained also
in \cite{SWZ} by methods of q-differential calculus.) 

Let us denote by \ $\cU_1$\ the unknown yet dual algebra of $\cA_1$, 
and by \ $\tA,B,C,\tD$\ the four generators of \ $\cU_1$. We would
like as in \cite{Dod} to define the pairing \ $\lg Z, f\rg$, \ $Z=\tA,B,C,\tD$,
\ $f$\ is from (\ref{eq:PBWbb}), as the classical tangent vector at the
identity: 
\eqn{duda} 
\lg \ Z \ , \ f \ \rg \ \equiv \ 
\ve \left( 
\frac{\partial f}{\partial y} \right) \ , \end{equation} 
however, here this would work only for the pairs: 
$\ (Z,y) \ =\ (\tA,\ta),\ (B,b),\ (\tD,\td)\ $, but not for $\
(C,c)\ $. The reason is that classically some of the relations in
(\ref{nrel}) are constraints and we have to differentiate
internally with respect to the manifold described by these 
constraints. In particular, if a constraint is given by setting
$\ g =0\ $, where $\ g\ $ is some function of $\ \ta,b,c,\td\ $, then any
differentiation $\ \cD\ $ should respect: 
\eqn{intd} \left(\cD\,g\,f\right)_{g=0}\ =\ 0 \ ,\end{equation} 
where $\ f\ $ is any polynomial function of $\ \ta,b,c,\td$. 
Thus, we are lead to define: 
\begin{equation}
 \left\langle C, f \rg 
\ \equiv\ \ve \left( E\ \frac{\partial}{\partial c} f\right)
\end{equation}
where: 
\eqna{eudb}
&&E\ =\ {\hat E}(- \tilde a, \frac{\partial}{\partial \tilde a}) \
, \\ 
&& {\hat E}(x,y)\ \equiv\ \sum_{k=0}^{\infty} \frac{x^k y^k}{k!} 
\eena

{}From the above definitions we get: 
\eqna{dudb}
&&\lg \ \tA \ , \ f \ \rg \ =\ 
\ve\left(\frac{\partial f}{\partial \ta}\right) \ =\ 
 \d_{n0} \begin{cases}k \qquad &{\rm for}\ \ f=b^n \ta^k \\ 
0 \qquad &{\rm for}\ \ f=b^n \td^\ell \\
0 \qquad &{\rm for}\ \ f=b^n c 
\end{cases} \\ 
&&\lg \ B \ , \ f \ \rg \ =\ 
\ve\left(\frac{\partial f}{\partial b}\right) \ =\ 
\d_{n1} \begin{cases}1 \qquad &{\rm for}\ \ f=b^n \ta^k \\ 
0 \qquad &{\rm for}\ \ f=b^n \td^\ell \\
0 \qquad &{\rm for}\ \ f=b^n c 
\end{cases} \\ 
&&\lg \ C \ , \ f \ \rg \ =\ 
\ve\left(E\, \frac{\partial f}{\partial c}\right) \ =\ 
\d_{n 0} \begin{cases}0 \qquad &{\rm for}\ \ f=b^n \ta^k \\ 
0 \qquad &{\rm for}\ \ f=b^n \td^\ell \\
1 \qquad &{\rm for}\ \ f=b^n c 
\end{cases} \\ 
&&\lg \ \tD \ , \ f \ \rg \ =\ 
\ve\left(\frac{\partial f}{\partial \td}\right) \ =\ 
\d_{\ell 1} \d_{n0} \begin{cases}
0 \qquad &{\rm for}\ \ f=b^n \ta^k \\ 
1 \qquad &{\rm for}\ \ f=b^n \td^\ell \\
0 \qquad &{\rm for}\ \ f=b^n c 
\end{cases} \qquad\\
&& \left\langle E, f \right\rangle \ =\ 
\begin{cases}1 &{\rm for}\ \ f=1_\cA\\ 
0 &{\rm otherwise} 
\end{cases} 
 \eena\nt 
We have included above also the auxiliary generator $\ E\ $
since it will appear in the coproduct relations (cf. below).
Note that if we have taken the definition (\ref{duda}) for $\
(C,c)\ $ the result in (\ref{dudb}) would superficially be the
same.

Now we can find the relations between the generators of
\ $\cA_1$. We have: 

\medskip 

\noindent 
{\bf Proposition 1:}\ The generators $\ \tA,B,C,\tD,E\ $
introduced above obey the following relations:
\eqna{cma} 
&& [\tilde D, C] \ =\ -2C \\
&& [B,C] \ =\ \tilde D \\
&& [B,C]_+ \ =\ \tD^2 \\
&& [\tilde D, B] \ =\ 2 B {\tilde D}^2 \\
&& [\tilde D, B]_+ \ =\ 0 \\
&& \tD^3 \ =\ \tD \\ 
&& C^2 \ =\ 0 \\ 
&& [\tA,B] \ =\ 0 \\ 
&& [\tA,C] \ =\ 0 \\ 
&& [\tA,\tD] \ =\ 0 \\ 
&& E Z \ =\ Z E \ =\ 0 \ , \quad Z=\tA,B,C,\tD \ . 
\eena 
\PR Using the assumed duality the above relations are
shown by calculating their pairings with the 
basis monomials \ $f$\ of $\cA_1\,$. In particular, we have:
\eqna{cmmu} 
&&\lg \ C \tD \ , \ f \ \rg \ =\ 
\begin{cases} 1 &{\rm for}\ \ f= c\\ 
0 &{\rm otherwise} 
\end{cases} \\
&&\lg \ \tD C\ , \ f \ \rg \ =\ 
\begin{cases} -1 &{\rm for}\ \ f= c\\ 
0 &{\rm otherwise} 
\end{cases}\\
&&\lg \ B C \ , \ f \ \rg \ =\ \begin{cases} 
\ha &{\rm for}\ \ f= \ta\\ 
\ha &{\rm for}\ \ f= \td\\ 
0 &{\rm otherwise} 
\end{cases} \\
&&\lg \ C B\ , \ f \ \rg \ =\ \begin{cases}
\ha &{\rm for}\ \ f= \ta\\ 
-\ha &{\rm for}\ \ f= \td\\ 
0 &{\rm otherwise} 
\end{cases}\\ 
&&\lg \ B \tD \ , \ f \ \rg \ =\ 
\begin{cases} -1 &{\rm for}\ \ f=b\\ 
0 &{\rm otherwise} 
\end{cases} \\
&&\lg \ \tD B\ , \ f \ \rg \ =\ \lg \ B\tD^2\ , \ f \ \rg \ =\ 
\begin{cases}1 &{\rm for}\ \ f=b\\ 
0 &{\rm otherwise} 
\end{cases}\\
&&\lg \ \tD^2 \ , \ f \ \rg \ =\ 
\begin{cases} 1 &{\rm for}\ \ f=a\\ 
0 &{\rm otherwise} 
\end{cases} \\ 
&&\lg \ \tA B \ , \ f \ \rg \ =\ \lg \ B \tA \ , \ f \ \rg \ =\ 
\begin{cases} k+1 &{\rm for}\ \ f=b \ta^k\\ 
0 &{\rm otherwise} 
\end{cases} \\ 
&&\lg \ \tA C \ , \ f \ \rg \ =\ \lg \ C \tA \ , \ f \ \rg \ =\ 
\begin{cases} 1 &{\rm for}\ \ f=c\\ 
0 &{\rm otherwise} 
\end{cases} \\ 
&&\lg \ \tA \tD \ , \ f \ \rg \ =\ \lg \ \tD \tA \ , \ f \ \rg \ =\ 
\begin{cases} 1 &{\rm for}\ \ f=\td\\ 
0 &{\rm otherwise} 
\end{cases}
 \eena 
\nt The Proposition now follows by formulae (\ref{cmmu})
and the defining relations (\ref{dudb}).$\ \diamondsuit$

We note that the algebraic relations (\ref{cma}) for $\un$ do not depend on
the constant $h$ present in the relations (\ref{nrel}) of the dual algebra
$\an$. Later, we shall see that the established duality reduces
also the algebra $\an$ so that it also does not depend on $h$. 

\subsection{Coalgebra structure of the dual}

We turn now to the coalgebra structure of \ $\un$. We have: 

\medskip 

\nt
{\bf Proposition 2:}\ \ 
(i) \ The comultiplication in the algebra \ $\un$\ is given by: 
\eqna{coa} 
&& \delta(\tilde A) \ =\ \tilde A\otimes 1_\cU + 1_\cU \otimes \tilde A, \\
&& \delta(B) \ =\ B\otimes 1_\cU + 1_\cU\otimes B, \\
 && \delta(C) \ =\ C\otimes E + E\otimes C, \\
 && \delta(\tilde D) \ =\ \tilde D\otimes E + E\otimes \tilde D, \\
 && \delta(E) \ =\ E\otimes E \eena 
(ii) \ The co-unit relations in \ $\un$\ are given
by: 
\eqna{coug} 
&&\ve_\cU (Z) \ =\ 0 \ , \qquad Z \ =\ \tA, B, C, \tD \\ 
&&\ve_\cU (E) \ =\ 1 \eena 
\nt where we have included also the auxiliary operator $E$.\\ 
\PR 
(i) \ We use the duality property (\ref{dub}{a}), namely we have 
$$\lg \ Z \ , \ f_1\ f_2 \ \rg \ =\ \lg \ \d_\cU(Z) \ , \ f_1
\otimes f_2 \ \rg $$
for every generator $Z$ of $\un$ and for every \ $f_1, f_2\in
\an$. Then we calculate separately the LHS 
and RHS and comparing the results prove (\ref{coa}).\nl 
(ii) \ Formulae (\ref{coug}) follow from \ $\ve_\cU (Z) \ =\
\lg \ Z , 1_\cA \ 
\rg$, cf. (\ref{dub}{b}), and using the defining relations
(\ref{dudb}).$\ \diamondsuit$ 

There is no antipode for the bialgebra $\un$. Indeed, suppose
that there was such. Then by applying the Hopf algebra axiom
(\ref{haxi}) to the generator $E$ we would get:
$$ E\, \g(E)\ =\ 1_\cU $$ 
which would lead to contradiction after multiplication from the
left with $\ Z \ =\ \tA, B, C, \tD\ $ (we would get $0=Z$).

\subsection{Reduction of the bialgebra} 

We noticed that the algebraic relations (\ref{cma}) of $\un$ 
do not depend on the constant $h$ from relations (\ref{nrel}) of
$\an$. The coproduct relations (\ref{coa}) also do not depend on $h$. 
We now clarify the reason for this. First we note that $\an$ 
has the following two-sided ideals and coideals:
\eqna{ideals}
&& I \ =\ {\an} b \tilde d \oplus {\an} \td^2 \oplus {\an} b c \\
&& I_2 \ =\ {\an} \td^2 \oplus {\an} b c \\
&& I_1 \ =\ {\an} b c \eena 
\nt so that 
\eqn{subi} I_1 \subset I_2 \subset I \subset \an 
\end{equation} 

Furthermore the pairing of all these ideals with the dual algebra
$\un$ vanish, thus we can set them consistently equal to zero.
Thus, the basis of $\an$ is reduced to the following monomials:
\eqn{newb} b^n \ta^k\, ,\ n,k\in \ZZ_+\,,\quad \td,\quad c 
\end{equation} 
Actually, it were only these monomials that appeared in the proof
of the dual relations (\ref{cma}). The algebraic relations of the
reduced algebra become rather trivial:
\begin{eqnarray}\label{nnrel}
&& 
\tilde{a}c \ =\ c\tilde{a} \ =\
\tilde{d}c \ =\ c\tilde{d} \ = \ \tilde{a}\tilde{d} \ =\
\tilde{d}\tilde{a} \ =\ cb \ =\ bc \ =\ 
\tilde{d}b \ =\ b\tilde{d} \ =\ 0 , \nn\\ 
&& c^2 \ =\ 0, \ \ \ \ \tilde{a}b \ =\ b\tilde{a} \ ,
\eea 
while the coalgebra relations remain unchanged and nontrivial. 
It is remarkable that the dual algebra has much richer structure
both in the algebraic and coalgebraic sectors.

\subsection{Consistency with the FRT approach} 

For the application of the FRT approach to duality we need the 
$4\times 4$ R-matrix which for the algebra $\an$ is given by
(\ref{rhu}). In the duality relations enter
actually the matrices $\ R_1^\pm\ $:
\eqna{rupm} 
&& R_1^+\ \equiv\ P\,R_1\,P\ =\ R_1 (-h)\ =
\ \Rtriang{-h}{h}{h_3}
 {1}{0}{h}
 {1}{-h}
 {1}
 \\ \nn\\ 
&& R_1^-\ \equiv\ R_1^{-1}\ =
\ \Rtriang{-h}{h}{-h_3-2h^2}
 {1}{0}{h}
 {1}{-h}
 {1}
 \eena \nt 
where $P$ is the permutation matrix:
\eqn{perm} 
P\ \equiv\ \left(
 \begin{array}{cccc}
 1 & 0 & 0 & 0 \cr
 0 & 0 & 1 & 0 \cr
 0 & 1 & 0 & 0 \cr
 0 & 0 & 0 & 1 \cr	 
 \end{array}
 \right)
\end{equation} 
 
 These R-matrices encode (part of) the duality between $\un$ and
$\an$ by formula (2.1) of \cite{FRT} taken for $k=1$ and written
in our setting: 
\eqn{frtd} \lg\ L^\pm \ ,\ T \ \rg \ =\ R_1^\pm \ ,\end{equation}
where $\ L^\pm\ $ are $2\times 2$ matrices whose elements are
functions of the generators of $\un$, $T$ is the $2\times 2$ matrix
formed by the generators of $\an$, c.f., (\ref{coop}). 
In order to make formula (\ref{frtd}) explicit we have to adopt some
convention on the indices. We choose to write it as: 
\eqn{frtdu} \lg\ L^\pm_{ik}\ ,\ T_{\ell j} \ \rg 
\ =\ (R_1^\pm)_{ijk\ell} \ ,\quad i,j,k,\ell =1,2\ ,\end{equation} 
where the enumeration of the R-matrices is done as in \cite{CFFS},
namely the rows are enumerated from top to bottom by the pairs
$\ (i,j)\ =\ (1,1), (1,2), (2,1), (2,2)$, and the columns 
are enumerated from left to right by the pairs $\ (k,\ell) \ =\ 
(1,1), (1,2), (2,1), (2,2)$. 

Using all this and rewriting the result in terms of the new basis
(\ref{nrel}) of $\an$ we have:
\eqn{loo} 
 \left\langle L^{\pm}_{11}\ ,\ 
 \left(
 \begin{array}{cc}
 \tilde a & b \cr c & \tilde d
 \end{array}
 \right)
 \right\rangle
 \ =\ 
 \left\langle L^{\pm}_{22}\ ,\ 
 \left(
 \begin{array}{cc}
 \tilde a & b \cr c & \tilde d
 \end{array}
 \right)
 \right\rangle
 \ =\ 
 \left(
 \begin{array}{cc}
 1 & - h \cr 0 & 0
 \end{array}
 \right)
\end{equation}
\eqn{lot} 
 \left\langle L^{\pm}_{12}\ ,\ 
 \left(
 \begin{array}{cc}
 \tilde a & b \cr c & \tilde d
 \end{array}
 \right)
 \right\rangle
 \ =\ 
 \left(
 \begin{array}{cc}
 h & h_\pm \cr 0 & 0
 \end{array}
 \right) \ ,
\end{equation}
where $\ h_+ \ =\ h_3\ $ and $\ h_- \ =\ -h_3 -2h^2$\ .
Note that the elements $\ L^{\pm}_{21}\ $ have zero products with
all generators so we can set them to zero. Next we calculate the
pairings with arbitrary elements of $\an$ for which we use the
fact that the coproducts of the $L^{\pm}_{jk}$ generators are
canonically given by \cite{FRT}:
\eqn{frtco} \d \left( L^{\pm}_{ik} \right) \ =\ \sum_{j=1}^2\ 
L^{\pm}_{ij}\ \otimes\ L^{\pm}_{jk} \ . \end{equation} 
Using this we obtain: 
\eqn{loab} 
 \left\langle L^{\pm}_{11} \ ,\ 
 b^n \tilde a^k
 \right\rangle
 \ =\ 
 \left\langle L^{\pm}_{22} \ ,\ 
 b^n \tilde a^k
 \right\rangle
 \ =\ (-h)^n
\end{equation}
\eqn{ltab}
 \left\langle L^{\pm}_{12} \ ,\ 
 b^n \tilde a^k
 \right\rangle
 \ =\ (-1)^n h^{n-1} ((k+n)h^2 - n (h_\pm + h^2))
\end{equation} 
All other pairings are zero. 

Computing the above pairings with the defining relations
(\ref{dudb}) we conclude that these $L$ operators are expressed
in terms of the generators of the dual algebra $\un$ as follows:
\eqna{lexp} 
 && L^{\pm}_{11} \ =\ L^{\pm}_{22} \ =\ e^{-h B} \\
 && L^{\pm}_{12} \ =\ ((h_\pm + h^2)B + h \tilde A) e^{-h B} 
\eena \nt
where expressions like $\ e^{\nu B}\ $ 
are defined as formal power series $\ e^{\nu B}\ =\ 1_\cU + 
\sum_{p\in\ZZ_+} \frac{\nu^p}{p!}\, B^p\ $.
Formulae (\ref{lexp}) are compatible with the coproducts 
(\ref{coa}{a,b}) of the generators $\ \tA,B\ $. However, 
as we see this approach does not say anything about the
generators $\ C,\tD\ $.

\section{Exotic bialgebras: case 2} 
\setcounter{equation}{0}

\subsection{Bialgebra relations}

In this Section we consider the bialgebra, denoted here by
$\ \cA_2\ $, which is obtained by applying the basic relations
(\ref{rtt}) for the case when $\ R\ =\ R_2\ $: 
\eqn{rhr} 
R_2\ =\ \Rtriang{h_1}{h_2}{h_3}
 {1}{0}{h_2}
 {1}{h_1}
 {1} 
 \quad , \qquad h_1 + h_2 \neq 0 
\end{equation} 
This $R$-matrix together with the condition on the parameters is the
second of the special cases (mentioned in the Introduction) of
the $R$-matrix denoted by $\ R_{H2,3}\ $ in
\cite{Hietarinta}. Its algebraic relations thus obtained 
are given by formulae (5.9) of \cite{AACDM}, namely: 
\begin{eqnarray}
 \label{eq:RTTh23a}
 &&c^2 \ =\ 0, \qquad ca \ =\ ac\ =\ 0, 
 \qquad dc \ =\ cd\ =\ 0, \nonumber\\
 && da \ =\ ad, \qquad cb \ =\ bc, \nonumber\\
 && a^2 \ =\ d^2\ =\ ad+bc, \nonumber\\ 
 && ab \ =\ bd \ = \ ba + (h_1-h_2) bc, 
 \qquad
 db \ =\ bd + (h_2-h_1) bc\ . 
\end{eqnarray} 
Note that the constant $h_3$ does not enter the above relations. 

The coalgebra relations are the same as for $\an$. Also the
demonstration that this bialgebra is not a Hopf algebra is
done as for $\,\cA_1\,$. 
The PBW basis in this case is:
\begin{equation}
 \label{eq:PBWa}
b^n a^k \,, \quad a^\ell d\,, \quad c\,, \qquad 
n,k\in \ZZ_+ \,, \quad \ell=0,1.
\end{equation}

Also in this case we make the change of basis (\ref{subs}) to
obtain:
\begin{eqnarray}\label{rrel}
&&c^2 \ =\ 0, \ \ \ \ \tilde{a}c \ =\ c\tilde{a} \ =\
\tilde{d}c \ =\ c\tilde{d} \ = \ \tilde{a}\tilde{d} \ =\
\tilde{d}\tilde{a} \ =\ 0 \ ,\nn\\ 
&& \ta b = b \ta \ ,\nn\\ 
&& b c \ =\ c b \ =\ 2\td^2 \ ,\qquad \td^3\ =\ 0 
\nn\\ 
&& \td b \ =\ - b \td \ =\ (h_1-h_2) \td^2 
\eea 
The PBW basis becomes:
\begin{equation}
 \label{eq:PBWaa}
b^n \ta^k \,, \quad \td^\ell\,, \quad c\,, \qquad 
n,k\in \ZZ_+ \,, \quad \ell=1,2 \,. 
\end{equation}
Thus, this bialgebra looks 'smaller' than $\ \an\ $ - compare with
(\ref{eq:PBWaa}). It has also a
smaller structure of two-sided ideals and coideals:
\eqna{idealz}
&& I_2 \ =\ {\ar} \td^2 \oplus {\ar} b c \\
&& I_1 \ =\ {\ar} b c \eena 
\nt so that 
\eqn{subii} I_1 \subset I_2 \subset \ar \end{equation} 
- compare with (\ref{ideals},\ref{subi}).

\subsection{Algebra and coalgebra structure of the dual}

In view of the similarities between the algebras $\an$ and $\ar$ 
it is natural do use the same generators $\ \tA,B,C,\tD,E\ $ for
the dual $\ur$. It is not surprising that we get the same
algebraic and coalgebraic relations. We have: 

\medskip 

\noindent 
{\bf Proposition 3:}\ The generators $\ \tA,B,C,\tD,E\ $ of the
bialgebra $\ \ur\ $ obey the same algebraic and coalgebraic
relations as for the algebra $\ \un\ $ given in Propositions 1
and 2.\\ \PR 
The proof is based on the fact that the bialgebras $\ \an\ $ and
$\ \ar\ $ differ in the relations involving the (co)ideals $\
I_k\ $ which have no bearing on the relations of $\ \un\ $. Thus,
we need only to show that all bilinears built from the generators
$\ \tA,B,C,\tD,E\ $ have zero pairings with the ideals $\ I_k\ $,
 (cf. (\ref{idealz},\ref{subii})), 
which is easy to demonstrate.\dia 

As a corollary also here the basis and algebraic relations of $\
\ar\ $ reduce to (\ref{newb}), (\ref{nnrel}). 

Thus, we have the following important conclusion:

\medskip 

\noindent 
{\bf Proposition 4:} The bialgebras $\ \an\ $ and $\ \ar\ $ considered as 
bialgebras in duality with the bialgebras $\ \un\ $, $\ \ur\ $,
respectively, coincide. 

We recall that the notion of duality we use does not coincide with the
FRT definition of duality. The latter is more stringent as we
shall see in the next subsection.

\subsection{Consistency with the FRT approach} 

The $4\times 4$ R-matrix needed for the FRT approach is given in
(\ref{rhr}). The matrices $\ R_2^\pm\ $ entering the duality
relations are: 
\eqna{rrpm} 
&& R_2^+\ \equiv\ P\,R_2\,P\ =\ 
\Rtriang{h_2}{h_1}{h_3}
 {1}{0}{h_1}
 {1}{h_2}
 {1}
 \\ \nn\\ 
&& R_2^-\ \equiv\ R_2^{-1}\ =
\ \Rtriang{-h_1}{-h_2}{2h_1h_2-h_3}
 {1}{0}{-h_2}
 {1}{-h_1}
 {1}
 \eena 
 
 Using the above and relations (\ref{frtdu}) (with $R_1\to R_2$)
we obtain: 
\begin{equation}
 \left\langle L^{+}_{11} , 
 \left(
 \begin{array}{cc}
 \tilde a & b \cr c & \tilde d
 \end{array}
 \right)
 \right\rangle
 \ =\ 
 \left\langle L^{+}_{22} , 
 \left(
 \begin{array}{cc}
 \tilde a & b \cr c & \tilde d
 \end{array}
 \right)
 \right\rangle
 \ =\ 
 \left(
 \begin{array}{cc}
 1 & h_2 \cr 0 & 0
 \end{array}
 \right)
\end{equation}
\begin{equation}
 \left\langle L^{+}_{12} , 
 \left(
 \begin{array}{cc}
 \tilde a & b \cr c & \tilde d
 \end{array}
 \right)
 \right\rangle
 \ =\ 
 \left(
 \begin{array}{cc}
 h_1 & h_3 \cr 0 & 0
 \end{array}
 \right)
\end{equation}
 \begin{equation}
 \left\langle L^{-}_{11} , 
 \left(
 \begin{array}{cc}
 \tilde a & b \cr c & \tilde d
 \end{array}
 \right)
 \right\rangle
 \ =\ 
 \left\langle L^{-}_{22} , 
 \left(
 \begin{array}{cc}
 \tilde a & b \cr c & \tilde d
 \end{array}
 \right)
 \right\rangle
 \ =\ 
 \left(
 \begin{array}{cc}
 1 & -h_1 \cr 0 & 0
 \end{array}
 \right)
\end{equation}
\begin{equation}
 \left\langle L^{-}_{12} , 
 \left(
 \begin{array}{cc}
 \tilde a & b \cr c & \tilde d
 \end{array}
 \right)
 \right\rangle
 \ =\ 
 \left(
 \begin{array}{cc}
 - h_2 & -h_3+2h_1 h_2 \cr 0 & 0
 \end{array}
 \right)
\end{equation}

Iterating this we obtain: 
\begin{equation}
 \left\langle L^{+}_{11} , 
 b^n \tilde a^k
 \right\rangle
 \ =\ 
 \left\langle L^{+}_{22} , 
 b^n \tilde a^k
 \right\rangle
 \ =\ h_2^n
\end{equation}
\begin{equation}
 \left\langle L^{+}_{12} , 
 b^n \tilde a^k
 \right\rangle
 \ =\ h_2^{n-1} ((k+n)h_1 h_2 + n (h_3-h_1 h_2))
\end{equation}
\begin{equation}
 \left\langle L^{-}_{11} , 
 b^n \tilde a^k
 \right\rangle
 \ =\ 
 \left\langle L^{-}_{22} , 
 b^n \tilde a^k
 \right\rangle
 \ =\ (-h_1)^n
\end{equation}
\begin{equation}
 \left\langle L^{-}_{12} , 
 b^n \tilde a^k
 \right\rangle
 \ =\ (-h_1)^{n-1} ((k+n)h_1 h_2 + n (-h_3+h_1 h_2))
\end{equation} 

{}From the above follow: 
\eqna{rexpr} 
 && L^{+}_{11} \ =\ L^{+}_{22} \ =\ e^{h_2 B} \\
 && L^{+}_{12} \ =\ ((h_3-h_1 h_2)B + h_1 \tilde A) e^{h_2 B} \\
 && L^{-}_{11} \ =\ L^{+}_{22} \ =\ e^{-h_1 B} \\
 && L^{-}_{12} \ =\ ((-h_3+h_1 h_2)B - h_2 \tilde A) e^{-h_1 B}
\eena 
This is compatible with the coproducts for the operators $\tA,B$.

Thus, we see that the $L$ operators in this case are different 
from those of $\ \un\ $, cf. (\ref{lexp}). Thus, the FRT approach
is more stringent than the notion of duality we use since it
distinguishes the two pairs of bialgebras. 
However, this difference is not as drastic as the difference
between the algebraic relations (\ref{nrel}), (\ref{rrel}) 
of $\ \an\ $, $\ \ar\ $, respectively, since (\ref{lexp}) is just
a special case of (\ref{rexpr}) obtained for $\ h_1= -h_2 =h\ $. 
On the other hand the FRT approach is incomplete in the cases at
hand since it gives info only about part of the generators,
namely, $\tA$ and $B$, and says nothing about the generators $C,\tD$.

\section{Exotic bialgebras: case 3} 
\setcounter{equation}{0}

\subsection{Bialgebra relations}

In this Section we consider the bialgebra which we denote here by
$\ \cA_3\ $. It is obtained by applying the basic relations
(\ref{rtt}) for the case when $\ R\ =\ R_3\ $: 
\eqn{rht} 
R_3\ =\ \Rtriang{0}{0}{1}
 {-1}{0}{0}
 {-1}{0}
 {1}
\end{equation} 
This R-matrix is denoted by $\ R_{S0,2}\ $ in \cite{Hietarinta}.
The algebraic relations of $\ \at\ $ are given by formulae (5.13)
of \cite{AACDM}, namely: 
\begin{eqnarray}
 \label{eq:RTTs02a}
 &&c^2\ =\ 0, \qquad ca \ =\ ac\ =\ 0, 
 \qquad dc \ =\ cd\ =\ 0, \nonumber\\
 && da\ =\ ad, \qquad cb \ =\ bc, \qquad a^2 \ =\ d^2 \nonumber\\ 
 && ab +ba\ =\ 0 , \qquad
 db + bd\ =\ 0 
\end{eqnarray} 
The coalgebra relations and the demonstration that this bialgebra
is not a Hopf algebra are as for $\an\,,\at$. 

Also in this case we make the change of basis (\ref{subs}) to
obtain:
\begin{eqnarray}\label{trel}
&&c^2 \ =\ 0, \ \ \ \ \tilde{a}c \ =\ c\tilde{a} \ =\
\tilde{d}c \ =\ c\tilde{d} \ = \ \tilde{a}\tilde{d} \ =\
\tilde{d}\tilde{a} \ =\ 0 ,\ \ \ \ cb \ =\ bc, \nn\\ 
&& \tilde{a}b \ +\ b\tilde{a}\ =\ 0 , \ \ \ \ \tilde{d}b \ +\
b\tilde{d} \ =\ 0 \ . \eea 

The algebra $\at$ has the same PBW bases (\ref{eq:PBWb}) and 
(\ref{eq:PBWbb}) as the algebra $\an$. It has also the same
(co)ideals as $\ \an\ $ (cf. (\ref{ideals},\ref{subi})).

\subsection{Algebra and coalgebra structure of the dual} 
 
In view of the similarities between the algebras $\an$ and $\at$
it is natural do use the same generators $\ \tA,B,C,\tD,E\ $ for
the dual $\ut$. It is not surprising that we get the same algebraic relations
between generators $\ \tA,B,C,\tD,E\ $. However, unlike the
bialgebras $\an$, $\ar$ the coalgebraic
relations and the relation with the FRT formalism here are 
different and it is even necessary to introduce two new 
auxilliary operators $\ F_\pm\ $ defined as:
\eqn{fiax} 
 \left\langle\ F_\pm, f\ \right\rangle 
 \ \equiv\ \ve \left({\hat E}(\pm 1, \frac{\partial}{\partial
\tilde d}) \, f\right) \ =\ 
\ve \left( \exp (\pm \frac{\partial}{\partial \tilde d}) 
\, f\right) \ .
\end{equation}
Explicitly we have:
\eqna{fiux} 
&& \left\langle\ F_+, f\ \right\rangle \ =\ 
\begin{cases} 1 &{\rm for}\ \ f=\td^\ell\\ 
1 &{\rm for}\ \ f=1_\cA \\ 
0 &{\rm otherwise} 
\end{cases} \\
&& \left\langle\ F_-, f\ \right\rangle \ =\ 
\begin{cases} (-1)^\ell &{\rm for}\ \ f=\td^\ell\\ 
1 &{\rm for}\ \ f=1_\cA \\ 
0 &{\rm otherwise} 
\end{cases} \eena 

We have for the algebraic and coalgebraic structure of $\ \ut\ $:

\medskip 

\noindent 
{\bf Proposition 5:}\ The generators $\ \tA,B,C,\tD,E,F_\pm\ $ 
 obey the following algebraic relations:
\eqna{cmat} 
&& [\tilde D, C] \ =\ -2C \\
&& [B,C] \ =\ \tilde D \\
&& [B,C]_+ \ =\ \tD^2 \\
&& [\tilde D, B] \ =\ 2 B {\tilde D}^2 \\
&& [\tilde D, B]_+ \ =\ 0 \\
&& \tD^3 \ =\ \tD \\ 
&& C^2 \ =\ 0 \\ 
&& [\tA,B] \ =\ 0 \\ 
&& [\tA,C] \ =\ 0 \\ 
&& [\tA,\tD] \ =\ 0 \\ 
&& E Z \ =\ Z E \ =\ 0 \ , \quad Z=\tA,B,C.\tD \\ 
&& F_+^2 \ =\ F_-^2 \ =\ 1_\cU \\ 
&& [F_+, F_-] \ =\ 0 \\ 
&& [\tA,F_\pm] \ =\ 0 \\ 
&& B F_\pm \pm F_\mp B \ =\ 0 \\ 
&& [C,F_\pm]_+ \ =\ 0 \\ 
&& [\tD,F_\pm] \ =\ 0 \\ 
&& E F_\pm \ =\ F_\pm E \ =\ E 
\eena 
\PR The relations between the generators $\ \tA,B,C,\tD,E\ $ are
the same as for the algebra $\un$ and they are proved in the same
way using (\ref{cmmu}). The only difference may have being in
(\ref{cmmu}{h}), since $\ b\ $ and $\ \ta\ $ anticommute, but
what is essential is that $\ \delta(b {\tilde a}^k)\ $ is
symmetric when paired with $\ \tA B\ $ and $\ B \tA\ $ (the
asymmetric terms involve $\ b c\ $ and $\ \td^2\ $ and give no
contribution). Verifying the relations involving $\ F_\pm\ $ is a
straightforward calculation.$\ \diamondsuit$ 

\medskip 

\noindent 
{\bf Proposition 6:}\ \ 
(i) \ The comultiplication in the algebra \ $\ut$\ is given by: 
\eqna{coat} 
&& \delta(\tilde A) \ =\ \tilde A\otimes 1_\cU + 1_\cU\otimes \tilde A, \\
&& \delta(B) \ =\ B\otimes 1_\cU + F_+F_-\otimes B, \\
&& \delta(C) \ =\ C\otimes E + E\otimes C, \\
&& \delta(\tilde D) \ =\ \tilde D\otimes E + E\otimes \tilde D, \\
&& \delta(E) \ =\ E\otimes E \\ 
&& \delta(F_\pm) \ =\ F_\pm \otimes F_\pm 
 \eena 
(ii) \ The co-unit relations in \ $\ut$\ are given
by: 
\eqna{cout} 
&&\ve_\cU (Z) \ =\ 0 \ , \qquad Z \ =\ \tA, B, C, \tD \\ 
&&\ve_\cU (Z) \ =\ 1 \ , \qquad Z \ =\ E, F_\pm \eena 
\PR 
The Proof is by the same methods as that of Proposition 2.\dia 

There is no antipode for the bialgebra $\ut$ - this is proved
exactly as for $\un$. 

As in the case of $\ \un\ \lra\ \an\ $ (and $\ \ur\ \lra\ \ar\ $)
duality one may reduce
the basis of $\ \at\ $ from the $\ \ut\ \lra\ \at\ $ duality, but
only with the ideal $\ I_1\ = \at\ bc $ (since $\ \td^2\ $ is not
annihilated by $\ F_\pm\ $). Thus, the basis of $\at$ is reduced
to the following monomials: 
\eqn{newbb} 
b^n \ta^k\,, \quad b^n \td^\ell \,, \quad c\,, \qquad 
n,k\in \ZZ_+ \,, \quad \ell\in\NN .
\end{equation} 
The algebraic relations of the reduced algebra become:
\begin{eqnarray}\label{ttrel}
&&c^2 \ =\ 0, \ \ \ \ \tilde{a}c \ =\ c\tilde{a} \ =\
\tilde{d}c \ =\ c\tilde{d} \ = \ \tilde{a}\tilde{d} \ =\
\tilde{d}\tilde{a} \ =\ cb \ =\ bc\ =\ 0 \ , \nn\\ 
&& \tilde{a}b \ +\ b\tilde{a}\ =\ 0 , \ \ \ \ \tilde{d}b \ +\
b\tilde{d} \ =\ 0 \ . \eea

\subsection{Consistency with the FRT approach} 

The $4\times 4$ R-matrix needed for the FRT approach is given in
(\ref{rht}). The matrices $\ R_3^\pm\ $ entering the duality
relations are: 
\eqna{rtpm} 
&& R_3^+\ \equiv\ P\,R_3\,P\ =\ R_3 \\ \nn\\ 
&& R_3^-\ \equiv\ R_3^{-1}\ =
\ \Rtriang{0}{0}{-1}
 {-1}{0}{0}
 {-1}{0}
 {1}
 \eena 
 
Using the above and relations (\ref{frtdu}) (with $R_1\to R_3$)
we obtain: 
\begin{equation}
 \left\langle L^{\pm}_{11} , 
 \left(
 \begin{array}{cc}
 \tilde a & b \cr c & \tilde d
 \end{array}
 \right)
 \right\rangle
 \ =\ 
 \left(
 \begin{array}{cc}
 0 & 0 \cr 0 & 1
 \end{array}
 \right)
\end{equation}
\begin{equation}
 \left\langle L^{\pm}_{22} , 
 \left(
 \begin{array}{cc}
 \tilde a & b \cr c & \tilde d
 \end{array}
 \right)
 \right\rangle
 \ =\ 
 \left(
 \begin{array}{cc}
 0 & 0 \cr 0 & -1
 \end{array}
 \right)
\end{equation}
\begin{equation}
 \left\langle L^{\pm}_{12} , 
 \left(
 \begin{array}{cc}
 \tilde a & b \cr c & \tilde d
 \end{array}
 \right)
 \right\rangle
 \ =\ 
 \left(
 \begin{array}{cc}
 0 & \pm 1 \cr 0 & 0
 \end{array}
 \right)
\end{equation}

Iterating these relations for arbitrary elements of the basis of
$\at$ we can show that the $\ L\ $ generators are given in terms
of some of the other generators in the following way: 
\eqn{mexp} 
 L^{\pm}_{11} \ =\ F_+\;, \qquad L^{\pm}_{22} \ =\ F_-\;, \qquad
 L^{\pm}_{12} \ =\ \pm B F_- 
\end{equation}
Formulae (\ref{mexp}) are compatible with the coproducts in 
(\ref{coat}) of the generators $\ B, F_\pm\ $. However, 
as we see this approach does not say anything about the basic 
generators $\ \tA, C,\tD\ $.

\section{Higher order R-matrix relations and quantum planes} 
\setcounter{equation}{0}

In order to address the question of the quantum planes 
corresponding to the exotic bialgebras we have to know the 
relations which the R-matrices fulfil. As we know the R-matrices 
producing deformations of the $GL(2)$ and $GL(1|1)$ fulfil second
order relations. However, in the cases at hand we have higher 
order relations. 

We start with the R-matrix $\ R_{H2,3}\ $ of \cite{Hietarinta}: 
\eqn{rhuu} 
R\ =\ \Rtriang{h_1}{h_2}{h_3}
 {1}{0}{h_2}
 {1}{h_1}
 {1}
\end{equation} \nt 
We need actually the singly permuted R-matrix:
\eqn{rhtt} 
\hR\ \equiv\ P\,R\ =\ 
\left(
 \begin{array}{cccc}
 1 & h_1 & h_2 & h_3 \cr
 0 & 0 & 1 & h_1 \cr
 0 & 1 & 0 & h_2 \cr
 0 & 0 & 0 & 1 
 \end{array}
 \right)
\end{equation} 

Explicit calculation shows now that we have:
\eqna{hwrm} 
&& ( \hR - \id )\, (\hR + \id)\ = \ 0 \ , \qquad h_1 = -h_2 =
h\ , \ h_3 = -h^2 \ , \\ 
&& ( \hR - \id )^2\, (\hR + \id)\ = \ 0 \ , \qquad h_1 = -h_2 =
h\ , \ h_3 \neq -h^2 \ , \qquad \hR\ =\ P\,R_1 \qquad \qquad \\ 
&& ( \hR - \id )^3\, (\hR+ \id)\ = \ 0 \ , \qquad h_1 + h_2 \neq
0 \ , \qquad \hR\ =\ P\,R_2 \ ,
\eena \nt 
where \ \id\ is the $4\times 4$ unit matrix. 
Thus the minimal polynomials are:
\eqn{mpol} {\rm pol}(\hR) = \begin{cases} 
( \hR - \id )\, (\hR + \id)  \qquad &{\rm for}\ h_1 = -h_2 =
h\ , \ h_3 = -h^2 \ , \\ 
( \hR - \id )^2\, (\hR + \id)\qquad &{\rm for}\ h_1 = -h_2 =
h\ , \ h_3 \neq -h^2 \\
 ( \hR - \id )^3\, (\hR+ \id) \qquad &{\rm for}\ h_1 + h_2 \neq
0 \\ 
\end{cases} \end{equation} 

\medskip 

\nt 
{\it Remark:}\ We recall that (\ref{hwrm}{a}) is 
the Jordanian subcase which produces the $GL_{h,h}(2)$ deformation
of $GL(2)$. Thus, the three subcases of Hietarinta's 
$R$-matrix $\ R_{H2,3}\ $ are distiguished not only and not so
much by the algebras they produce but intrinsically by their  
minimal polynomials.~\dia 

\medskip 

To derive the corresponding quantum planes we shall apply the
formalism of \cite{WZ} (cf. also \cite{Schwenk}). The commutaion relations 
between the coordinates $\ z^i\ $ and differentials $\ \z^i\ $, 
($i=1,2$), are given as follows:
\eqn{qpl} z^i z^j \ =\ {\cP}_{ijk\ell}\, z^k z^\ell \end{equation} 
\eqn{dfa} \zeta^i \zeta^j \ =\ - {\cQ}_{ijk\ell}\, \zeta^k \zeta^\ell 
 \end{equation} 
\eqn{cda} z^i \zeta^j \ =\ {\cQ}_{ijk\ell}\, \zeta^k z^\ell 
\end{equation} 
where the operators $\ \cP,\, \cQ\ $ are functions of $\ \hR\ $
and must satisfy:
\eqn{orth} (\cP - \id)\,(\cQ + \id) \ = \ 0\ . \end{equation} 
In the well studied deformations of $\ GL(2)\ $ there 
are quadratic minimal polynomials and there are only two
choices for the operators $\ \cP,\, \cQ\ $, cf. e.g.,
(\ref{hwrm}{a}). Here we have more choices. In particular, for
the case (\ref{hwrm}{b}) we have four 
choices: 
\eqn{choicea}
\left( \cP - \id\ ,\ \cQ + \id \right) = 
\begin{cases} 
\left( \hR - \id\ ,\ \hR^2 - \id \right) \\
\left( \hR + \id \ ,\ (\hR - \id)^2 \right) \\ 
\left( \hR^2 - \id\ ,\ \hR - \id \right) \\
\left( (\hR - \id)^2\ ,\ \hR + \id \right) 
\end{cases}
\end{equation} 
while in the case (\ref{hwrm}{c}) we have six choices: 
\eqn{choiceb}
\left( \cP - \id\ ,\ \cQ + \id \right) = 
\begin{cases} 
\left( \hR - \id\ ,\ (\hR^2 - \id)(\hR - \id)
 \right) \\
\left( \hR + \id \ ,\ (\hR - \id)^3 \right) \\ 
\left( \hR^2 - \id\ ,\ (\hR - \id)^2 \right) \\
\left( (\hR - \id)^2\ ,\ \hR^2 - \id \right)\\ 
\left( (\hR^2 - \id)(\hR-\id)\ ,\ \hR - \id \right) \\
\left( (\hR - \id)^3\ ,\ \hR + \id \right) 
\end{cases}
\end{equation} 

Our choice will be the last possibility of both (\ref{choicea}), 
\ref{choiceb}), i.e., we shall use $\ \cP - \id = (\hR - \id)^a\
$ with $a=2,3$, respectively, and $\ \cQ = \hR\ $ in all cases.
With this choices and 
denoting $\ (x,y) = (z^1,z^2)\ $ we obtain from (\ref{qpl}) 
\eqn{qpc} x y - y x = h y^2\ , \qquad h_1 = - h_2=h \ , \qquad
\cP - \id = (\hR - \id)^2\ , 
\end{equation} 
or
\eqn{qpb} 
x y - y x = \ha (h_1-h_2) y^2\ , \qquad h_1 \neq - h_2 \ , \qquad
\cP - \id = (\hR - \id)^3\ . 
\end{equation} 
We note that the quantum planes corresponding to the bialgebras
$\ \an\ $ and $\ \ar\ $ are not essentially different.
Furthermore the quantum plane (\ref{qpc}) is the same as for the
Jordanian subcase if we choose $\ \cP - \id = \hR - \id\ $.

Denoting $\ (\xi,\eta) = (\zeta^1,\zeta^2)\ $ we obtain from
(\ref{dfa}) with $\ \cQ=\hR\ $: 
\eqna{dfb}
&& \xi^2 + \frac{h_1 - h_2}{2}\ \xi \eta \ =\ 0 \\ 
&& \eta^2 \ =\ 0 \\ 
&& \xi \eta \ =\ - \eta \xi 
\eena \nt
Of course, for $\ \hR=P\,R_1\ $ (\ref{dfb}{a}) simplifies to 
\eqn{dfbb} \xi^2 + h\, \xi \eta \ =\ 0 \ ,
\end{equation} 
which is valid also for the Jordanian subcase. 

Finally, for the coordinates-differentials relations we obtain
from (\ref{cda}) with $\ \cQ=\hR\ $ again for all subcases: 
\eqna{cdb}
&& x \xi \ =\ \xi x + h_1 \xi y + 
h_2 \eta x + h_3 \eta y \\ 
&& x \eta \ =\ \eta x + h_1 \eta y \\ 
&& y \xi \ =\ \xi y + h_2 \eta y \\ 
&& y \eta \ =\ \eta y \eena

Finally we derive the quantum plane relations for the case of the $\ R_3\
$ matrix. It is easy to see that (\ref{hwrm}{b}) holds also in
this case, i.e., for 
\eqn{rhtz} 
\hR_3\ \equiv\ P\,R_3\ =\ 
\left(
 \begin{array}{cccc}
 1 & 0 & 0 & 1 \cr
 0 & 0 & -1 & 0 \cr
 0 & -1 & 0 & 0 \cr
 0 & 0 & 0 & 1 
 \end{array}
 \right)
\end{equation}

Using (\ref{qpl},\ref{dfa},\ref{cda}) with $\ \cP -\id =
(\hR_3-\id)^2\ $, $\ \cQ = \hR_3\ $, we obtain, respectively:
\eqn{qpd} 
x y \ =\ -y x \end{equation} 
\eqna{dfc}
&& \xi^2 \ =\ 0 \\ 
&& \eta^2 \ =\ 0 \\ 
&& \xi \eta \ =\ \eta \xi \eena 
\eqna{cdc}
&& x \xi \ =\ \xi x + \eta y \\ 
&& x \eta \ =\ - \eta x \\ 
&& y \xi \ =\ - \xi y \\ 
&& y \eta \ =\ \eta y \eena 

Finally, we note that a check of consistency of this formalism is
to implement Manin's approach to quantum planes \cite{Manin}. 
Namely, one takes quantum matrix $\ T\ $, cf. (\ref{coop})
as transformation matrix of the two-dimensional quantum planes.
This means, that if we define: 
\eqn{qpll} 
{z'}^i \ =\ T_{ij}\, z^j \ , \qquad {\zeta'}^i \ =\ T_{ij}\, \zeta^j \ , 
\end{equation} 
then $\ (x',y')\ =\ ({z'}^1,{z'}^2)\ $ and 
$\ (\xi',\eta') = ({\zeta'}^1,{\zeta'}^2)\ $ 
should satisfy the same relations as $\ (x,y)\ $ and $\
(\xi,\eta)\ $. The latter statement may be used to recover 
the algebraic relations of the bialgebras. Namely, suppose, 
that relations (\ref{qpc}),(\ref{dfb}{b,c}),(\ref{dfbb}),(\ref{cdb}), 
or relations (\ref{qpb}),(\ref{dfb}),(\ref{cdb}), 
or relations (\ref{qpd}),(\ref{dfc}),(\ref{cdc}), 
hold for both 
$\ (x,y)\ $ and $\ (\xi,\eta)\ $, 
$\ (x',y')\ $ and $\ (\xi',\eta')\ $; then substitute the 
expressions for $\ (x',y')\ $ and $\ (\xi',\eta')\ $ in the 
these relations, under the assumption that $\ a,b,c,d\ $ commute
with $\ (x,y)\ $ and $\ (\xi,\eta)\ $; then the coefficients of
the independent bilinears that may be built from $\ (x,y)\ $ and
$\ (\xi,\eta)\ $, will reproduce the algebraic relations of the
bialgebras $\ \an\ $, $\ \ar\ $, $\ \at\ $, respectively.

\section{Conclusions and outlook} 
\setcounter{equation}{0}

In this paper we have found the bialgebras which are in duality
with the three exotic matrix bialgebras (obtained in 
\cite{AACDM}) which are ~{\it not}~ deformations of the classical algebra 
of functions over the group $\ GL(2)$ or the supergroup
$GL(1|1)$. 

These bialgebras are rather degenerate and on their example we
discover several hitherto unknown phenomena. To illustrate this
we comment in more detail on the first two cases (considered in
sections 3 and 4). The starting point are the $4\times 4$
$R$-matrices $R_1$ and $R_2$, cf. (\ref{rhu}), (\ref{rhr}). On
the one hand $R_1$ is a special case of $R_2$ obtained for
\eqn{sbss} h_1=-h_2 =h \end{equation} 
On the other hand the algebraic relations obtained in 
\cite{AACDM} by applying the FRT formalism are different and
those for $\an$ (\ref{nrel}) cannot be obtained from those for $\ar$ 
(\ref{rrel}) by using (\ref{sbss}). Another peculiarity is that
the parameter $h_3$ does not appear neither in (\ref{nrel}), nor
in (\ref{rrel}). However, the dependence on the parameters turns
out to be redundant in both cases, after we find the bialgebras
they are in duality with, namely, $\un$ and $\ur$, respectively.
This is so since the parameter 
dependence enters only through the (co)ideals of the bialgebras,
which (co)ideals have zero pairings with $\un$ and $\ur$,
respectively. Thus, these (co)ideals can be neglected (or we can
pass to factor-algebras). As a result the reduced
bialgebras coincide and the same holds for their bialgebras in duality.
The two algebras differ, though only in a limited sense, if we 
use also the $\ L-T\ $ duality of FRT \cite{FRT}, cf.
(\ref{frtd}). The limited sense being that the $L$ operators have
the same parameter dependence as the matrices $R_1$ and $R_2$,
respectively, and in the same way the $L$ operators for 
$\un$ can be obtained from the $L$ operators of $\ur$ by 
(\ref{sbss}). 

Thus, the FRT formalism is more stringent since it preserves the
initial parameter dependence. On the other hand, it reproduces
only two of the basic generators of $\un$ and $\ur$, namely,
those, that we denote by $\tA,B$, and gives no information on the
other two basic generators $C,\tD$. This insufficiency of the FRT
formalism is similar to the one observed for the Jordanian
deformations, cf., e.g., \cite{Ohn}.

To conclude, the real difference between the two cases is
exhibited by only the minimal polynomials of the permuted
$R$-matrices. The significance of this will be revealed
fully in the representation theory of the exotic bialgebras and
their duals which we intend to develop next.

\paragraph{Acknowledgments:} This work was supported in part by
the CNRS-BAS\\ France/Bulgaria agreement number 6608. 

\np

\end{document}